\documentclass[12pt]{amsart}
\usepackage{amssymb}

\usepackage{color}
\newtheorem{theorem}{Theorem}[section]
\newtheorem{lemma}[theorem]{Lemma}

\newtheorem{proposition}[theorem]{Proposition}

\newtheorem{lem-def}[theorem]{Lemma-Definition}

\newcommand{\hooklongrightarrow}{\lhook\joinrel\longrightarrow}
\renewenvironment{proof}{{\bfseries Proof.}}{\qed}
\newcommand{\B}{\mathbb B}

\newcommand{\Z}{\mathbb Z}

\newcommand{\D}{\mathbb D}

\newcommand{\T}{\widetilde{\mathbb V}}
\newcommand{\V}{\mathbb V}

\def\op{\operatorname}

\def\al{\alpha}

\def\ap{\mathbf{A}}
\def\ars#1{\renewcommand\arraystretch{#1}}

\def\aut{\op{Aut}}

\def\bad{B(a,\delta)}
\def\bb{{\mathcal B}}
\def\bcut{\B_{\cut}}
\def\be{\beta}
\def\bg{{\mathbb{B}_\g}}
\def\bi{\B_\infty}

\def\bqc{\B_{\op{qcut}}}

\def\bul{B^{\mbox{\tiny $\bullet$}}}

\def\cc{\mathcal{C}}

\def\cut{{\op{cut}}}
\def\cuts{\op{Cuts}(\g)}
\def\inn{\op{in}}

\def\dcut{\D_\cut}

\def\defn{\nn{\bf Definition. }}

\def\dfd{D(f,\dta)}
\def\dg{\D_\g}

\def\di{\D_\infty}

\def\dqc{\D_{\op{qcut}}}

\def\diso{\lower.4ex\hbox{$\downarrow$}\raise.4ex\hbox{\mbox{\scriptsize
$\wr$}}}

\def\dta{\delta}
\def\dul{D^{\mbox{\tiny $\bullet$}}}

\def\e{\medskip}

\def\ep{\epsilon}

\def\g{\Gamma}
\def\ga{\gamma}

\def\gd{\g(\dta)}

\def\gg{\mathcal{G}}

\def\gm{\g_\mu}

\def\go{\g_\omega}

\def\gv{\Gamma_v}

\def\hk{\hookrightarrow}

\def\hra{\hooklongrightarrow}

\def\imp{\ \Longrightarrow\ }

\def\inii{\op{Init}(\g)}

\def\inn{\op{in}}

\def\irr{\op{Irr}}
\def\ism{\lower.3ex\hbox{\ars{.08}$\begin{array}{c}\,\to\\\mbox{\tiny $\sim\,$}\end{array}$}}
\def\iso{\lower.3ex\hbox{\ars{.08}$\begin{array}{c}\lra\\\mbox{\tiny $\sim\,$}\end{array}$}}

\def\kb{\overline{K}}
\def\vt{\vb_\tau}
\def\kh{K^h}
\def\khx{K^h[x]}

\def\kp{\op{KP}}

\def\kpm{\op{KP}(\mu)}

\def\kx{K[x]}

\def\La{\Lambda}
\def\lc{\op{lc}}

\def\lg{l\raise.6ex\hbox to.2em{\hss.\hss}l}

\def\lra{\,\longrightarrow\,}

\def\lx{\operatorname{lex}}

\def\nn{\noindent}
\def\nnn{\mathcal{N}}

\def\om{\omega}
\def\omi{\omega_{-\infty}}

\def\orb{\hbox to  .3em{$\backslash$}\backslash}
\def\ord{\op{ord}}
\def\mub{\overline{\mu}}
\def\p{\mathfrak{p}}

\def\qc{\op{Qcuts}(\g)}

\def\kx{K[x]}

\def\sg{\sigma}

\def\sii{\ \Longleftrightarrow\ }

\def\supp{\op{supp}}

\def\vb{\bar{v}}
\def\vc{\V^{\op{com}}}

\def\vi{\V_\infty}
\def\vf{\V_{\op{fin}}}

\def\vrt{\V_{\op{rt}}}
\def\vt{\V_{\op{vt}}}

\newcounter{cs}
\stepcounter{cs}
\newcommand{\casos}{\begin{itemize}}
\newcommand{\fcasos}{\end{itemize}\setcounter{cs}{1}}

\newfont{\tit}{cmr12 scaled \magstep3}

\makeatletter \@namedef{subjclassname@2010}{  \textup{2010} Mathematics Subject Classification} 

\title[Geometric parametrization of valuations]{Geometric parametrization of valuations on a polynomial ring}
\author{Enric Nart}
\address{Departament de Matem\`{a}tiques,         Universitat Aut\`{o}noma de Barcelona,         Edifici C, E-08193 Bellaterra, Barcelona, Catalonia}
\email{nart@mat.uab.cat}
\author{Josnei Novacoski}
\address{Departamento de Matem\'{a}tica,         Universidade Federal de S\~ao Carlos, Rod. Washington Luís, 235, 13565--905, S\~ao Carlos -SP, Brazil}
\email{josnei@ufscar.br}
\thanks{Partially supported by grant PID2020-116542GB-I00  funded by the Spanish MCIN/AEI. 
During the realization of this project the second author was supported by two grants from Funda\c{c}\~ao de Amparo \`a Pesquisa do Estado de S\~ao Paulo (process numbers 2017/17835-9 and 2021/11246-7).}

\begin{document}

\begin{abstract}
We extend and prove a conjecture of Bengu\c{s}-Lasnier on the parametri\-zation of valuations on a polynomial ring by certain spaces of diskoids. 
\end{abstract}

\maketitle

\section*{Introduction}

Let $(K,v)$ be a valued field with value group $\gv=vK$.
Let  $\kh$ be a henselization of $(K,v)$, determined by the choice of an extension $\vb$ of $v$ to a fixed algebraic closure $\kb$ of $K$. The value group $\g:=\g_{\vb}$  of $\vb$ is the divisible closure of $\gv$.
Let us denote the decomposition group of $\vb$ by
$$
G=G_{\vb}=\left\{\sg\in \aut(\kb/K)\mid \vb\circ\sg=\vb\right\}=\aut(\kb/\kh).
$$

Let $\V$ be the set of equivalence classes of valuations on $K(x)$ whose restriction to $K$ is equivalent to $v$. In \cite[Section 4]{Rig}, it is shown that $\V$ has a natural structure  of a tree; that is, a partially ordered set whose  intervals are totally ordered.

Let $\vi\subseteq \V$ be the subset of classes of {\bf valuation-algebraic} valuations, which can be constructed only after some limit process. This set $\vi$ is included in the set of leaves (maximal nodes) of our tree $\V$. Let  $\vf\subseteq{\V}$ be the subset of classes of {\bf valuation-transcendental} valuations, which can be constructed from suitable  monomial valuations  by a finite number of augmentations. We have a decomposition
\[
\V=\vi\sqcup \vf=\vi\sqcup \vrt\sqcup\vt,
\]
where $\vrt$, $\vt$ are the subsets of classes of {\bf residue-transcendental} and {\bf value-transcendental} valuations, respectively. 

In a recent paper \cite{Andrei}, Bengu\c{s}-Lasnier conjectures the existence of a geometric space of diskoids parametrizing the subtree $\vrt$ of residue-transcendental valuations. This is proven in \cite{Andrei} in the henselian case, and in the rank-one case. 

The aim of this paper is to generalize and prove this conjecture for arbitrary valuations in $\V$, with no assumptions on the valued field $(K,v)$.
 
Following  Bengu\c{s}-Lasnier's insight, we construct a space $\D$ of certain gene\-ra\-lized diskoids, equipped with a natural partial ordering, and define a mapping:
\begin{equation}\label{maineq}
\D\lra\V,\qquad D\longmapsto \mu_D,
\end{equation}
which is an isomorphism of posets.

Also, we find concrete descriptions of subspaces $\di,\, \dg,\, \dcut\subseteq\D$ such that the above isomorphism induces isomorphisms  
\[
\di\iso \vi,\qquad \dg\iso \vrt,\qquad \dcut\iso \vt.
\]

The space $\dg$ contains all diskoids  centred at monic irreducible polynomials in $\kh[x]$ and whose radii  belong to $\g$. The space $\dcut$ contains more general diskoids admitting cuts in $\g$ as radii.

The space $\di$  is the quotient set of the space of all nests of diskoids in $\dg$  having an empty intersection, under the equivalence relation which identifies two nests if they are mutually cofinal.

The outline of the paper is as follows. Section \ref{secBackground} reviews some basic facts on commensurable extensions of valuations and valuative trees. In Section \ref{secKb}, we para\-me\-trize $\V$ when the field $K$ is algebraically closed. As shown by Alexandru-Popescu-Zaharescu \cite{APZ0,APZ}, the set  $\vrt$ is parametrized by ultrametric closed balls in $\kb$ with radii in $\g$. The description of $\vi$ in  terms of nests of closed balls with empty intersection  can be found in \cite{APZ2,Kuhl,Vaq3}. For the parametrization of $\vt$ we use balls in $\kb$ admitting cuts in $\g$ as radii.  

In Section \ref{secProof11} we prove our main theorem (Theorem \ref{main}). 
Clearly, we can parametrize valuations on $\kh(x)$ just by considering the $G$-orbits of the geometric spaces of balls (or nests of balls) parametrizing valuations on $\kb(x)$. Here,  $G$ is the decomposition group defined above. By the rigidity theorem \cite[Theorem B]{Rig}, we get a parametrization of  valuations on $K(x)$ as well. 
Finally, it is easy to identify $G$-orbits of balls with diskoids centred at irreducible polynomials in $\khx$.

Section \ref{secSupport} is devoted to extend these results to parametrize equivalence classes of valuations on the polynomial ring $\kx$, including those with nontrivial support.
This is achieved by admitting diskoids with an infinite radius.  

Finally, in Section \ref{secApprT}, we compare our geometric parametrization with that obtained by using Kuhlmann's  {\bf approximation types} \cite{KuhlAT}. For an algebraically closed field, Kuhlmann's parametrization is an isomorphism too, and  our mapping in (\ref{maineq}) is essentially equal to the inverse of the approximation types mapping. 
  
\section{Background on valuative trees}\label{secBackground}

\subsection{Commensurable extensions of valuations}\label{subsecComm}
Let $\vc$ be the set of all extensions of $v$ to $K(x)$ taking values in a fixed divisible closure $\g$ of $\gv$. We consider in $\vc$ the following partial ordering:
\[
\mu\le\nu\ \sii\ \mu(f)\le\nu(f)\ \mbox{ for all }f\in \kx.
\]
As  usual, we write $\mu<\nu$ when $\mu\le\nu$ and $\mu\ne\nu$. With this structure of a poset, $\vc$ becomes a tree; that is, all intervals in $\vc$ are totally ordered. 

A valuation $\mu$ on $K(x)$ extending $v$  is {\bf commensurable} (over $v$) if $\gm/\gv$ is a torsion group. 

Clearly, all valuations in $\vc$ are commensurable. Conversely, every commensurable extension of $v$ to $K(x)$ is equivalent to a unique valuation  in $\vc$. Indeed, for any given extension $\iota\colon\gv\hk\La$ of ordered abelian groups such that $\La/\gv$ is a torsion group, there is a unique embedding $\La\hk\g$ such that the composition $\gv\hk\La\hk\g$ is the canonical embedding of $\gv$ in $\g$. 
In particular, two valuations in $\vc$ are equivalent if and only if they coincide.

Since valuation-algebraic and residue-trans\-cen\-dental valuations are commensura\-ble, 
the sets of equivalence classes of these valuations  can be identified with subsets of $\vc$, and they fill up the whole space:
\[
\vc=\vi\sqcup \vrt.
\]

A valuation $\mu\in\vc$ is residue-transcendental if and only 
the residual extension $K(x)\mu/Kv$ is transcendental.
This is equivalent to $\kpm\ne\emptyset$, where $\kpm$ is the set of all Mac Lane--Vaqui\'e (MLV) key polynomials for $\mu$ \cite[Theorem 4.4]{KP}.

Therefore, in order to fulfill our aim, we must parametrize {\bf all} valuations in $\vc$ and, moreover, parametrize the set $\vt$ of  {\bf equivalence classes} of value-transcen\-den\-tal valuations. 

The latter task is much more subtle. To start with, the fact that there is a canonical partial ordering on $\vt$ is not obvious and will be discussed in the next section.

\subsection{Partial ordering on $\V$}\label{subsecPosetV}
Let $\mu$ be an arbitrary extension of $v$ to $K(x)$. By \cite[Theorem 1.5]{Kuhl}, the extension $\gv\hk\gm$ is {\bf small}. That is, if $\gv\subseteq\Delta\subseteq\gm$ is the relative divisible closure of $\gv$ in $\gm$, the quotient $\gm/\Delta$ is a cyclic group.

Two small extensions of $\gv$ are said to be $\gv$-equivalent if they are isomorphic as ordered abelian groups, by an isomorphism acting as the identity on $\gv$. 

In \cite[Section 4.2]{csme}, a minimal universal ordered abelian group $\g\subseteq \La$ is constructed, containing all small extensions of $\gv$ up to $\gv$-equivalence.
Hence, we have a natural partial ordering on the set of all $\La$-valued valuations on $K(x)$, which descends to a canonical partial ordering on $\V$ \cite[Section 4.1]{Rig}.

In practice, for any $\mu,\nu\in\V$, we may still think that
\[
\mu\le\nu\ \sii\ \mu(f)\le \nu(f)\quad \mbox{for all }f\in\kx,
\]   
where the inequality $\mu(f)\le\nu(f)$ must be properly understood. If $\mu\ne\nu$ (i.e. $\mu$ is not equivalent to $\nu$), then this inequality is equivalent to the existence of $\ga\in\g$ such that
$\mu(f)\le\ga\le \nu(f)$.

\section{Valuations on $\kb(x)$}\label{secKb}
As mentioned in the Introduction, our first step will be to parametrize valuations on $\kb(x)$ in terms of ultrametric balls in $\kb$. For any field $K\subseteq L\subseteq \kb$, let us write 
\[
\V(L)=\vi(L)\sqcup \vf(L)=\vi(L)\sqcup \vrt(L)\sqcup\vt(L),
\]
for the corresponding spaces of equivalence classes of valuations on $L(x)$.

\subsection{Residue-transcendental valuations}\label{subsecRTKb}  

The parametrization of $\vrt(\kb)$ is well-known \cite{APZ0,APZ,Kuhl,Vaq3}. Let us review the results omitting all proofs.

For each pair $(a,\dta)\in\kb\times\g$, consider the closed ball of center $a$ and radius $\dta$:
\[
\bad=\left\{c\in\kb\mid \vb(c-a)\ge\dta\right\}.
\] 
The criterion of coincidence of two balls is:
\[
\bad=B(b,\ep)\ \sii \ \vb(b-a)\ge\dta=\ep.
\]
In particular, any $b\in\bad$ is a center of the ball: $\bad=B(b,\dta)$.

Let $\bg$ be the set of all these ultrametric closed balls.

For every $B=B(a,\dta)\in\bg$, denote by $\om_B=\om_{a,\dta}$ the monomial valuation defined as follows in terms of $(x-a)$-expansions:
\[
\om_B\left(\sum\nolimits_{0\le n}a_n(x-a)^n\right) = \min\left\{\vb(a_n)+n\dta\mid0\le n\right\}.
\]
These valuations $\om_B$ admit the following reinterpretation. 

\begin{proposition}\cite[Lemma 6.2]{Andrei}\label{infRT}
	For all $B\in\bg$, we have
	\[
	\om_B(f)=\min\{\vb(f(c))\mid c\in B\}\quad\mbox{ for all }f\in\kb[x].
	\]
\end{proposition}

These extensions of $\vb$ to $\kb(x)$ are residue-transcendental. Actually, for all $c\in\kb$ we have \cite[Lemma 2.4]{Rig}:
\[
\ars{1.2}
\begin{array}{l}
c\in B\ \imp\ x-c\in\kp(\om_B),\\
c\not\in B\ \imp\ \inn_{\om_B}(x-c)\ \mbox{ is a unit in }\gg_{\om_B},
\end{array}
\] 
where $\gg_{\om_B}$ is the graded algebra associated to $\om_B$. On the other hand, all residue-transcendental valuations on $\kb(x)$ arise in this way \cite[Theorem 2.1]{APZ0}. 

Also, for all $B,C\in\bg$ we have
\[
\om_B=\om_C\ \sii\ B=C.
\]
Therefore, we obtain a bijective mapping
$$
\bg\lra\vrt(\kb),\qquad B\longmapsto \om_B. 
$$

Finally, consider the partial ordering  in $\bg$ determined by descending inclusion:
\[
B\le C\ \sii\ B\supseteq C.
\] 
It follows easily from Proposition \ref{infRT} that the mapping $B\mapsto\om_B$ strictly preserves the ordering:
\[
B<C\ \imp\ \om_B<\om_C.
\] 
Therefore, our bijective mapping is an isomorphism of posets.  

\begin{theorem}\label{mainRTKb}
	The mapping $\,\bg\to\vrt(\kb)$ determined by $B\mapsto \om_B$ is an isomorphism of posets.
\end{theorem}

\subsection{Valuation-algebraic valuations}\label{subsecVAKb}  

Let $\nnn_\infty$ be the set of all nests of closed balls of $\bg$, with empty intersection. Thus, an element in $\nnn_\infty$ is a family $\bb=(B_i)_{i\in I}$ of balls $B_i\in\bg$, parametrized by some totally ordered set $I$ of indices, such that
\[
i<j\ \sii \ B_i<B_j\ \sii\ B_i\supset B_j.
\]  
Moreover, $\bigcap_{i\in I}B_i=\emptyset$. In particular, the set $I$ contains no last element.

Denote for simplicity $\om_i=\om_{B_i}$ for all $i\in I$. The family $\Omega_\bb=(\om_i)_{i\in I}$ is a {\bf continuous  family} of valuations on $\kb(x)$ of degree one (every valuation admits a MLV key polynomial of degree one).
For the definition of continuous families and their limits see \cite{Vaq} or \cite[Section 4.1]{VT}.

The property $\bigcap_{i\in I}B_i=\emptyset$ is equivalent to the fact that $\Omega_{\bb}$ has a stable limit, which we denote as:
\[
\om_\bb=\lim\left(\Omega_\bb\right).
\]
This means that all polynomials $f\in\kb[x]$ are $\Omega_{\bb}$-stable; that is, there exists an index $i_0\in I$ (depending on $f$) such that
\[
\om_i(f)=\om_{i_0}(f) \ \mbox{ for all } i\ge i_0.
\]
Then, $\om_\bb(f)$ is defined to be this stable value. In other words:
\[
\om_\bb(f)=\max\{\om_i(f)\mid i\in I\}.
\]
These valuations $\om_\bb$ are valuation-algebraic \cite[Proposition 4.1]{VT}. Moreover, every valuation-algebraic valuation on $\kb(x)$ arises in this way \cite[Theorem 5.1]{APZ2}, \cite[Proposition 2.12]{Vaq3}.

On the set $\nnn_\infty$ we consider a natural equivalence relation. Two nests of closed balls $\bb=(B_i)_{i\in I}$, $\cc=(C_j)_{j\in J}$ are said to be equivalent if they are mutually cofinal. That is, for all $i\in I$ there exists $j\in J$ such that $B_i\le C_j$, and vice versa.
In this case, we write $\bb\sim \cc$.

Let $\bi=\nnn_\infty/\!\sim$ \,be the quotient set of $\nnn_\infty$ under this equivalence relation. Denote by $[\bb]\in\bi$ the class of any $\bb\in\nnn_\infty$.   

By Theorem \ref{mainRTKb}, the continuous families $\Omega_\bb$, $\Omega_\cc$ are mutually cofinal if and only if $\bb\sim\cc$. By \cite[Proposition 4.8]{VT}, we have
\[
\om_\bb=\om_\cc\ \sii\ \bb\sim \cc.
\]
The next result follows.

\begin{theorem}\label{mainVAKb}
	The mapping $\bi\to\vi(\kb)$, defined by $\,[\bb]\mapsto \om_\bb$ is bijective.
\end{theorem}

These valuations $\om_\bb\in\vi(\kb)$ follow the pattern indicated in Proposition \ref{infRT} too. 

\begin{proposition}\label{infVA}
For $\bb=(B_i)_{i\in I}\in\nnn_\infty$, let $B=\bigcup_{i\in I}B_i$. Then,
	\[
	\om_\bb(f)=\min\{\vb(f(c))\mid c\in B\}\quad\mbox{ for all }f\in\kb[x].
	\]
\end{proposition}

This follows easily from Proposition \ref{infRT} and the definition of $\om_\bb$.

\subsection{Cuts in $\g$}\label{subsecCuts}  

The equivalence classes of value-transcen\-den\-tal valuations on $\kb(x)$ can be parametrized  by balls admitting cuts in the group $\g$ as radii.

Let us recall some basic properties of cuts.
For all $\ga\in \g$, we denote 
\[
\g_{< \ga}=\{\al\in \g\mid \al< \ga\}\subset \g_{\le \ga}=\{\al\in \g\mid \al\le \ga\}.
\] 

 For $S,T\subseteq \g$ and $\ga\in \g$, the  expressions
	$$\ga<S,\ \quad  \ga\le S,\  \quad S<T,\ \quad S\le T$$
	mean that the corresponding inequality holds for all $\al\in S$ and all $\be\in T$.

An {\bf initial segment} of $\g$ is a subset $S\subseteq \g$ such that $\g_{\le \ga}\subseteq S$ for all $\ga\in S$.

On the set $\inii$ of initial segments of $\g$ we consider the ordering determined by ascending inclusion. We obtain a totally ordered set with a minimal and a maximal element:
$\emptyset=\min(\inii)$, $\g=\max(\inii)$.

A {\bf cut} in $\g$ is a pair $\dta=(\dta^L,\dta^R)$ of subsets of $\g$ such that $$\dta^L< \dta^R\quad\mbox{ and }\quad \dta^L\cup \dta^R=\g.$$ Clearly, $\dta^L$ is an initial segment of $\g$. 
Let $\cuts$ be the set of all cuts in $\g$. We have an obvious bijection
\[
\inii\lra\cuts,\qquad S  \longmapsto (S,\g\setminus S).
\]
We consider on $\cuts$ the total ordering induced by $\inii$  through this bijection.
In particular, $\cuts$ admits a minimal and a maximal element
\[
-\infty:=(\emptyset,\g)=\min\left(\cuts\right),	\qquad \infty^-:=(\g,\emptyset)=\max\left(\cuts\right),
\]  which are called {\bf improper cuts}.

The notation $\infty^-$ is motivated by the fact that this cut  is the immediate predecessor of $\infty$ in the totally ordered set $\cuts\infty$.\e
 
\defn
Every $\ga\in\g$ determines two {\bf principal cuts}:
\[
\ga^-=\left(\g_{<\ga},\,\g_{\ge\ga}\right)
,\qquad \ga^+=\left(\g_{\le\ga},\,\g_{>\ga}\right).
 \] \vskip0.1cm

For every cut $\dta=(\dta^L,\dta^R)\in\cuts$, consider a formal symbol $x_\dta$ and build up the abelian group 
$\gd= x_\dta\Z\oplus\g$.
There is a unique ordering on $\gd$ which is compatible with the group structure and satisfies
$\dta^L<x_\dta<\dta^R$.
Namely, 
$$
mx_\dta+b\le nx_\dta+a\sii (m-n)x_\dta\le a-b\sii (m-n)\dta^L\le a-b.
$$
The latter inequality means that  $(m-n)\ga\le a-b$, for all $\ga\in\dta^L$.

Clearly, the extension $\g\hk\gd$ is incommensurable. However, it is a small extension of ordered abelian groups because $\gd/\g$ is a cyclic group.  

The groups determined by the improper cuts have a specially simple description. Indeed, the following maps are isomorphisms of ordered abelian groups:
\begin{equation}\label{lex}
\ars{1.4}
\begin{array}{l}
\g(-\infty)\lra\left(\Z\times \g\right)_{\op{lex}},\qquad mx_{-\infty}+\ga\ \mapsto\ (-m,\ga),\\
\g(\infty^-)\lra\left(\Z\times \g\right)_{\op{lex}},\qquad\; mx_{\infty^-}+\ga\ \mapsto\ (m,\ga).
\end{array}
\end{equation}

\subsection{Value-transcendental valuations}\label{subsecValueTKb}  
For every pair $(a,\dta)\in\kb\times\cuts$, consider the ball of center $a$ and radius $\dta$, defined as:
\[
\bad=\left\{c\in\kb\mid \vb(c-a)>\dta^L\right\}\subseteq\kb.
\] 

For instance, the balls having a principal cut as radius are the standard closed and open balls with radius in $\g$. Indeed, for all $\ga\in\g$ we have:
\begin{equation}\label{NotMerge}
\ars{1.3}
\begin{array}{l}
B(a,\ga^-)=\left\{c\in\kb\mid \vb(c-a)\ge\ga\right\}=B(a,\ga),\\
B(a,\ga^+)=\left\{c\in\kb\mid \vb(c-a)>\ga\right\}=B^{\op{o}}(a,\ga).
\end{array}
\end{equation}

Also, the improper cuts determine very special balls. For all $a\in\kb$, we have
\begin{equation}\label{ImproperBalls}
B(a,-\infty)=\kb,\qquad B(a,\infty^-)=\{a\}.
\end{equation}

\begin{lemma}\label{coincidence}
	Let $a,b\in\kb$ and $\dta,\ep\in\cuts$. Then,
	\[
	\bad=B(b,\ep)\ \sii \ \vb(b-a)>\dta^L\ \mbox{ and }\ \dta=\ep.
	\]
\end{lemma}

\begin{proof}
	Suppose $\bad=B(b,\ep)$. Since $b\in \bad$, we have $\vb(b-a)>\dta^L$.
	
	If $\dta<\ep$, then $\dta^L\subsetneq \ep^L$ and there exists $\ga\in\g$ such that $\dta^L<\ga\le\ep^L$. If $c\in\kb$ has value $\vb(c)=\ga$, then we have $c\in \bad$ and $c\not\in B(b,\ep)$, contradicting our assumption. The inequality  $\ep<\dta$ leads to a completely analogous contradiction. Hence, $\dta=\ep$.
	The converse implication is obvious. 
\end{proof}\e

In particular, any $b\in\bad$ is a center of the ball: $\bad=B(b,\dta)$.

Let $\bcut$ be the set of all these ultrametric balls centred at elements in $\kb$ and having a cut as radius.

For every $B=B(a,\dta)\in\bcut$, denote by $\om_B=\om_{a,\dta}$ the monomial valuation defined as follows in terms of $(x-a)$-expansions:
\begin{equation}\label{defom}
\om_B\left(\sum\nolimits_{0\le n}a_n(x-a)^n\right) = \min\left\{\vb(a_n)+nx_\dta\mid0\le n\right\}.
\end{equation}
Note that $\om_B$ is value-transcen\-den\-tal with  value group $\g_{\om_B}=\gd$.

For all $c\in\kb$ we have:
\begin{equation}\label{x-c}
\ars{1.2}
\begin{array}{l}
c\in B\ \imp\ \om_B(x-c)=x_\dta,\\
c\not\in B\ \imp\ \om_B(x-c)=\vb(c-a)\in\g.
\end{array}
\end{equation}

As in the classical case, these valuations $\om_B$ are uniquely determined by their defining balls.

\begin{lemma}\label{B=om}
	For all $B,C\in\bcut$ the following conditions are equivalent.
	\begin{enumerate}
		\item[(a)] $B=C$. 
		\item[(b)] $\om_B=\om_C$.
		\item[(c)] $\om_B\sim\om_C$.
	\end{enumerate}
\end{lemma}

\begin{proof}
	Obviously, (a)$\Rightarrow$(b)$\Rightarrow$(c). Let us show that (c)$\Rightarrow$(a).
	
	Suppose that $B=\bad$, $C=B(b,\ep)$ and $\om_B\sim\om_C$. This means that  there is an isomorphism of ordered groups $\iota\colon \g(\ep) \ism\gd$ fitting into a commutative diagram
	
	$$
	\ars{1.4}
	\begin{array}{c}
	\g(\ep)\ \stackrel{\iota}\lra\ \gd\\
	\quad\ \mbox{\scriptsize$\om_C$}\nwarrow\quad \nearrow\mbox{\scriptsize$\om_B$}\quad\\
	K(x)^*
	\end{array}
	$$
	
	Since $\om_B$ and $\om_C$ are extensions of $\vb$, the isomorphism $\iota$ acts as the identity on $\g$. 
	
	For all $c\in C$, (\ref{x-c}) shows that $\om_C(x-c)=x_\ep$. Since the diagram commutes, 
	\[
	\om_B(x-c)=\iota\left(\om_C(x-c)\right)=\iota\left(x_\ep\right)\not\in\g.
	\] By (\ref{x-c}), we deduce that $c\in B$. This shows that $C\subseteq B$, and we deduce that $B=C$ by the symmetry of the argument.
\end{proof}\e

On the other hand, all value-transcendental valuations on $\kb(x)$ arise in this way. 

\begin{lemma}\label{ontoVT}
	Every value-transcendental valuation on $K(x)$ is equivalent to $\om_B$ for some $B\in\bcut$.  
\end{lemma}

\begin{proof}
	Let $\om$ be a valuation-algebraic extension of $\vb$ to $\kb(x)$. By \cite[Theorem 1.5]{Kuhl}, the value group $\go$ of $\om$ is an extension of $\g$ as an ordered abelian group of the following form:  
	\[
	\g\hra\go=\al\Z\oplus\g, \qquad \ga\longmapsto (0,\ga),
	\] where $\al\in\go$ has no torsion over $\g$. By \cite[Theorem 3.11]{Kuhl}, there exists $a\in \kb$ such that $\om$ is the monomial valuation acting as usual on $(x-a)$-expansions:
	\[
	\om\left(\sum\nolimits_{0\le n}a_n(x-a)^n\right) = \min\left\{\vb(a_n)+n\al\mid0\le n\right\}.
	\] 
	Let $\dta\in\cuts$ be the cut in $\g$ determined by $\al$. That is, 
	\[
	\dta^L=\left\{\ga\in\g\mid \ga<\al\right\},\qquad \dta^R=\left\{\ga\in\g\mid \ga>\al\right\}.
	\]
	We may buid up an isomorphism of ordered abelian groups
	$\iota\colon \gd\to\go$, acting as the identity on $\g$ and mapping $x_\dta$ to $\al$.  If $B=\bad$, then we have clearly $\om=\iota\circ\om_B$. Hence, $\om\sim\om_B$.
\end{proof}\e

Therefore, we obtain a bijective mapping
$$
\bcut\lra\vt(\kb),\qquad B\longmapsto \om_B. 
$$

Finally, consider the partial ordering  in $\bcut$ determined by descending inclusion:
\[
B\le C\ \sii\ B\supseteq C.
\] 
It is easy to check that the mapping $B\mapsto\om_B$ strictly preserves the ordering. Therefore, out bijective mapping is an isomorphism of posets.  

\begin{theorem}\label{mainVTKb}
	The mapping $\,\bcut\to\vt(\kb)$ determined by $B\mapsto \om_B$ is an isomorphism of posets.
\end{theorem}

\subsubsection{Minimal and maximal elements in $\vt(\kb)$}\label{subsubsecMM}
The tree $\vt(\kb)$ has an absolute minimal element, which we denote by $\omi$. Indeed, by (\ref{ImproperBalls}) and Theorem \ref{mainVTKb}, for all $a,b\in\kb$ we have
\begin{equation}\label{omi}
\omi:=\om_{a,-\infty}=\om_{b,-\infty}\le \om\quad\mbox{ for all }\om\in\vt(\kb).
\end{equation}

We say that $\om_{-\infty}$ is the {\bf root node} of $\vt(\kb)$. By the isomorphism in (\ref{lex}), we may think that $\omi$ takes values in the group $\left(\Z\times\g\right)_{\lx}$ and
acts as follows on nonzero polynomials $f\in\kb[x]$:
\[
\omi(f)=\left(-\ord_x(f),\vb(\lc(f))\right),
\] 
where $\lc(f)$ is the leading coefficient of $f$.

On the other hand, by (\ref{ImproperBalls}) and Theorem \ref{mainVTKb}, for all $a\in\kb$ the valuation $\om_{a,\infty^-}$ is maximal in $\vt(\kb)$. Again, by using the isomorphism in (\ref{lex}), we may think that $\om_{a,\infty^-}$ takes values in the group $\left(\Z\times\g\right)_{\lx}$ and
acts as follows on nonzero polynomials $f\in\kb[x]$:
\[
\om_{a,\infty^-}(f)=\left(\ord_{x-a}(f),\vb(\inn(f))\right),
\] 
where $\inn(f)$ is the first nonzero coefficient in the $(x-a)$-expansion of $f$.

\subsection{Valuation-transcendental valuations}\label{subsecQC}  

The tree $\vf=\vrt\sqcup \vt$ whose nodes represent equivalence classes of valuation-transcendental extensions of $v$ to $K(x)$ has been described in \cite[Section 7]{VT}. It is a ``compact" set, in the sense that any totally ordered subset has an infimum and a supremum. In particular,  every two nodes have a greatest common lower node. Also, it has a unique root node and as many maximal nodes as monic irreducible polynomials in $\kh[x]$.

After Theorems \ref{mainRTKb} and \ref{mainVTKb}, a geometric space parametrizing $\vf(\kb)$ should be a kind of merge of the spaces of balls $\bg$ and $\bcut$.  However, this merge cannot be a simple union, because $\bg\subseteq\bcut$ as subsets of $\mathcal{P}(\kb)$, as shown in (\ref{NotMerge}). 

In other words, each ball in $\bg$ determines two different valuations, depending on our consideration of the radius being some $\ga\in\g$ or the corresponding cut $\ga^-\in\cuts$. Even if we considered the formal disjoint union $\bg\sqcup \bcut$ as our space, there would remain the problem of deciding what partial ordering on this formal union reflects the partial ordering on $\vf(\kb)$. 

A natural merging of these spaces is obtained by considering balls admitting quasi-cuts in $\g$ as radii. 

A {\bf quasi-cut} in $\g$ is a pair $\dta=\left(\dta^L,\dta^R\right)$ of subsets of $\g$ such that 
\[
\dta^L\le \dta^R\quad \mbox{ and }\quad \dta^L\cup \dta^R=\g.
\]
Thus, $\dta^L$ is an initial segment of $\g$ and $\dta^L\cap \dta^R$ consists of at most one element.

The set $\qc$ of all quasi-cuts in $\g$ admits a total ordering:
$$
\dta=\left(\dta^L,\dta^R\right)\le \ep=\left(\ep^L,\ep^R\right) \ \sii\ \dta^L\subseteq \ep^L\quad\mbox{and}\quad \dta^R\supseteq \ep^R.
$$

There is an embedding of ordered sets $\g\hk\qc$, which assigns to every $\ga\in\g$ the {\bf principal quasi-cut} $\left(\g_{\le \ga},\g_{\ge \ga}\right)$. We abuse of language and denote still by $\ga\in \qc$ the principal quasi-cut determined by $\ga$.

Clearly, the set $\cuts$ is embedded in $\qc$ and 
\[
\qc=\g\sqcup\cuts.
\]  
As far as the ordering is concerned, the comparison between the elements of $\g$ and $\cuts$ is clarified by the following inequalities:
\[
\ga^-<\ga<\ga^+\quad\mbox{ for all }\ga\in\g.
\] 
Also, $\ga^-$ is the immediate predecessor of $\ga$ in $\qc$, while  $\ga^+$ is the immediate successor of $\ga$.

For each pair $(a,\dta)\in\kb\times\qc$, consider the ``pointed" ball of center $a$ and radius $\dta$, defined as:
\[\bul(a,\dta)=\left(\bad,\dta\right),
\] 
where $\bad\subseteq \kb$ is the ball defined in Sections \ref{subsecRTKb} (for $\dta\in\g$) and \ref{subsecValueTKb}  (for $\dta\in\cuts$). We denote by $\bqc$ the set of all these pointed closed balls.

Thus, a pointed closed ball $\bul\in\bqc$ has two ingredients: 

$\bullet$ \ an {\bf underlying ball} $B\subseteq \kb$,   denoted $B=\op{under}(\bul)$.

$\bullet$ \ a {\bf radius} in $\qc$.\e

\noindent{\bf Definition. }We say that $\,\dta\in\qc$ \,is an \emph{essential radius} if $\dta^R$ contains no minimal element. Equivalently, \,$\dta\not\in\g$ and $\dta\ne\ga^-$  \,for all \,$\ga\in\g$.\e

A pointed ball with an essential radius is uniquely determined by its underlying ball in $\kb$. However, for all $\ga\in\g$, the pointed balls $\bul(a,\ga)$ and $\bul(a,\ga^-)$  have the same underlying ball.

For every  $\bul=\left(\bul(a,\dta),\dta\right)\in\bqc$, we denote $\om_{\bul}=\om_{a,\dta}$.

In the set $\bqc$, we consider the following partial ordering:
\[
\bul(a,\dta)\le \bul(b,\ep)\ \sii \bad\supseteq B(b,\ep) \ \mbox{ and }\ \dta\le\ep.
\]

Then, it is easy to deduce from Theorems \ref{mainRTKb} and \ref{mainVTKb} the following result.

\begin{theorem}\label{mainQC}
The mapping $\bqc\to\vf(\kb)$ determined by $\bul\mapsto \om_{\bul}$ is an isomorphism of posets.
\end{theorem}

The following result generalizes Proposition \ref{infRT} to arbitrary valuation-transcenden\-tal valuations. For a quasi-cut $\dta\in\g\subset\qc$, let us write $x_\dta:=\dta\in\g$.

\begin{proposition}\label{weneed}
	Let $\om=\om_{a,\dta}$ for some $a\in \kb$, $\dta\in\qc$. Take a non-zero
	$f=\sum\nolimits_{0\le n}a_n(x-a)^n\in\kb[x]$.
	Let $S$ be the set of indices $n$ such that 
	$\om(f)=\vb(a_n)+nx_\dta$.
	Denote $V=\left\{\vb(f(c))\mid c\in B(a,\dta)\right\}$. 
	Then, the following holds:
	\begin{enumerate}
		\item[(a)] If $0\in S$ or $\dta\in\g$, then  $\ \om(f)=\min\left(V\right)$.
		\item[(b)] If $0\not\in S$ and $\dta=\ga^-$ for some $\ga\in\g$, then  $\ \om(f)=\min\left(V\right)^-$.
		\item[(c)] If $0\not\in S$ and $\dta$ is an essential radius, then $\om(f)=\inf\left(V\right)$.
	\end{enumerate}
\end{proposition}

\nn{\bf Remark. }In (b), an equality of the form $\om(f)=\al^-$, for $\al\in\g$, means that the value $\om(f)\in\go$ strictly realizes the cut $\al^-$; that is, $\g_{<\al}<\om(f)<\g_{\ge\al}$.

In (c), the equality  $\om(f)=\inf\left(V\right)$ means that $\om(f)$ is the infimum of $V$ as a subset of $\go$.\e

\begin{proof}
	Let us write $B=B(a,\dta)$. For all $c\in B$ we have
\[
	\vb(f(c))\ge \min\{\vb(a_n(c-a)^n)\mid 0\le n\}\ge \min\{\vb(a_n)+nx_\dta\mid 0\le s\}=\om(f).
\]	
We deduce that  $\om(f)\le V$.

	If $0\in S$, then $\om(f)=\vb(a_0)=\vb(f(a))$. Since $a\in B$, we have $\om(f)=\min\left(V\right)$.
	
	If $\dta\in\g$, then (a) follows from Proposition \ref{infRT}. This ends the proof of (a).
		
	From now on, we assume that $0\not\in S$ and $\dta\not\in\g$. Then, $S$ is a one-element set, because for all $n\ne m$ we have
	$$
	\vb(a_n)+nx_\dta=\vb(a_m)+mx_\dta\ \imp x_\dta=(\vb(a_n)-\vb(a_m))/(m-n)\in\g,
	$$
against our assumption. Hence, $S=\{n_0\}$ with $n_0>0$, and
	\begin{equation}\label{n0}
	\om(f)=\vb(a_{n_0})+n_0x_\dta<\vb(a_{n})+nx_\dta,\quad \mbox{for all }n\ne n_0. 
	\end{equation}
	
Since $n_0>0$, we have $\om(f)\not\in\g$, so that $\om(f)<V$.\e

\nn{\bf Proof of (b). }Suppose $0\not\in S$ and $\dta=\ga^-$ for some $\ga\in\g$. By (\ref{NotMerge}), $B=B(a,\ga)$. By Propositon \ref{infRT}, $\min(V)$ exists and coincides with $\om_{a,\ga}(f)$. Let us show that 
\[
\min(V)=\om_{a,\ga}(f)=\vb(a_{n_0})+n_0\ga.
\]
Indeed, since $\ga>x_\dta=x_{\ga^-}$, we have  $n(\ga-x_\dta)>n_0(\ga-x_\dta)$ for all $n>n_0$. Hence, we deduce from (\ref{n0}) that $\vb(a_{n_0}+n_0\ga)<\vb(a_n)+n\ga$. 
On the other hand, for $n<n_0$, the inequality (\ref{n0}) implies 
\[
x_\dta<\dfrac{\vb(a_n)-\vb(a_{n_0})}{n_0-n}. 
\]
Hence, $\ga\le (\vb(a_n)-\vb(a_{n_0}))/(n_0-n)$, because $\ga$ is the minimal element in $\dta^R$. This leads to  
$\vb(a_{n_0})+n_0\ga\le\vb(a_n)+n\ga$. Hence, 
\[
\om_{a,\ga}(f)=\min\{\vb(a_n)+n\ga\}=\vb(a_{n_0})+n_0\ga.
\]

In order to show that $\om(f)$ strictly realizes the cut $\min(V)^-$ in $\go$, it suffices to show that $\g_{<\min(V)}<\om(f)$. Consider any  $\al\in\g$ such that $\om(f)<\al$, and  express it as $\al=\vb(a_0)+n_0\be$ for some $\be\in\g$. We have 
\[
\vb(a_0)+n_0x_{\ga^-}=\om(f)<\al=\vb(a_0)+n_0\be\ \imp\ x_{\ga^-}<\be\ \imp\ \ga\le\be. 
\]
Hence, $\al\ge\vb(a_0)+n_0\ga=\min(V)$. This ends the proof of (b).\e

\nn{\bf Proof of (c). }Suppose $0\not\in S$ and $\dta$ is an essential radius.
Write 
\[
\al:=\min\{\left(\vb(a_n)-\vb(a_{n_0})\right)/(n_0-n)\mid 0\le n< n_0\}\in\g.
\]

\noindent{\bf Claim. }Suppose that $\be\in\g$ satisfies $x_\dta<\be<\al$. Then, for any $c\in\kb$ such that $\vb(c-a)=\be$, we have $c\in B$ and $\om(f)<\vb(f(c))=\vb(a_{n_0})+n_0\be$.\e

Indeed,  $c\in B$ because $\vb(c-a)>x_\dta>\dta^L$.  Also, for all $n<n_0$, we have
\[
\vb(c-a)<\dfrac{\vb(a_n)-\vb(a_{n_0})}{n_0-n}\ \imp\ \vb(a_n(c-a)^n)<\vb(a_{n_0}(c-a)^{n_0}).
\]
Finally, for all $n>n_0$, from (\ref{n0}) and $x_\dta<\be$,  we  deduce
$ \vb(a_{n_0})+n_0\be<\vb(a_{n})+n\be$.
Hence, $\vb(a_n(c-a)^n)<\vb(a_{n_0}(c-a)^{n_0})$ for all $n>n_0$.
As a consequence,
\[
\vb(f(c))=\vb(a_{n_0}(c-a)^{n_0})=\vb(a_{n_0})+n_0\be>\vb(a_{n_0})+n_0x_\dta=\om(f).
\]
This ends the proof of the Claim.

In order to show that $\om(f)=\inf(V)$, we need only to check that for every $\xi\in\go$ such that $\om(f)<\xi$, there is some $\ep\in V$ such that $\om(f)<\ep<\xi$. In other words, there is some $c\in B$ such that $\om(f)<\vb(f(c))<\xi$. 

Since $\dta$ is an essential radius and $n_0>0$, the cut of $\g$ determined by $\om(f)=\vb(a_0)+n_0x_\delta$ is essential too. Hence, we may assume that $\xi\in\g$.

In order to facilitate some comparisons, let us write  $\xi=\vb(a_{n_0})+n_0\ga$ for some $\ga\in\g$. Then, the inequality $\om(f)<\xi$ is equivalent to $x_\dta<\ga$. Since $x_\dta<\al$ and $\dta^R$ has no minimal element, there exists $\be\in\g$ such that $x_\dta<\be<\min\{\ga,\al\}$. 

By the Claim, there exists $c\in B$ such that 
\[
\om(f)<\vb(f(c))=\vb(a_{n_0})+n_0\be<\vb(a_{n_0})+n_0\ga=\xi.
\]
This ends the proof of (c).
\end{proof}

\subsection{Global ordering on $\B$}\label{subsecGlobalOrd}

Let $\B=\bi\sqcup \bqc$. Since the intersection of a nest of closed balls in $\bi$ is empty, for any reasonable extension of the partial ordering of $\bqc$ to $\B$, the nests of closed balls in $\bi$ must be maximal elements.

Thus, we obtain a partial ordering on $\B$ just by defining
what balls in $\bqc$ lie below any given nest of balls.

Let $\bb=(B_i)_{i\in I}\in\nnn_\infty$ be a nest of closed balls, and $\bul\in\bqc$. Then, we define
\[
\bul\le \bb\ \sii\ \op{under}(\bul)\le B_i\ \mbox{ for some }i\in I.  
\] 
Clearly, this ordering is compatible with the equivalence $\sim$ of nests of closed balls. Thus, it determines a partial ordering on $\B$.

We deduce our main theorem for algebraically closed fields.

\begin{theorem}\label{mainB}
The mapping $\,\B\to\V(\kb)$ determined by $B\mapsto \om_B$ is an isomorphism of posets.
\end{theorem}
	
The node $\omi\in\vt(\kb)\subseteq\V(\kb)$ becomes a root node of $\V(\kb)$ too.

Also, for all $a\in\kb$, the nodes $\om_{a,\infty^-}\in\vt(\kb)\subseteq\V(\kb)$ become maximal nodes of $\V(\kb)$ too. We say that these $\om_{a,\infty^-}$ are {\bf finite leaves} of $\V(\kb)$.

The set of leaves of $\V(\kb)$ is
\[
\bi\sqcup\{\om_{a,\infty^-}\mid a\in\kb\}.
\]
We say that the leaves in $\bi$ are {\bf infinite leaves} of $\V(\kb)$.

Every interval between the root node and a finite leaf  is parametrized by $\qc$:
\[
[\omi,\om_{a,\infty^-}]=\left\{\om_{a,\dta}\mid \dta\in \qc\right\}.
\] 

\section{Descent of valuations from $\kb(x)$ to $K(x)$}\label{secProof11}

The descent of valuations from $\kb(x)$ to $\kh(x)$ is described in full detail in a recent paper by Vaqui\'e \cite[Section 3]{Vaq3}.
Since the automorphisms in the decomposition group $G$ leave $\vb$ invariant, they act on ultrametric balls as follows:
\[
\sg\left(\bad\right)=B(\sg(a),\dta),\quad \forall\,\sg\in G,\ \forall\,\dta\in\qc.
\]
Since this action preserves the radius $\dta$, it extends in an obvious way to an action on the set $\bqc$ of pointed balls. Also, $G$ acts on $\bi$ because its action on nests of closed balls is clearly compatible with the equivalence relation of mutual cofinality. 

Therefore, $G$ acts in a natural way on $\B$. Let $\B/G$ be the set of all $G$-orbits, and let us denote by $[B]_G$ the $G$-orbit of any $B\in\B$.

This $G$-action on $\B$ is reflected in the following action on $\V$:
\[
\om_{\sg(B)}=\om_B\circ\sg^{-1}\quad\mbox{ for all }\sg\in G.
\]

Hence, for all $B,C\in\B$ we have
\begin{align*}
\left(\om_B\right)_{\mid \kh(x)}=\left(\om_C\right)_{\mid \kh(x)}&\ \sii\ \om_C=\om_B\circ\sg\ \mbox{ for some }\sg\in G\\&\ \sii\ C = \sg^{-1}(B)\ \mbox{ for some }\sg\in G\\&\ \sii \ [B]_G=[C]_G.
\end{align*}

Therefore, we deduce from Theorem \ref{mainB} an isomorphism of posets
\[
\B/G\lra\V(\kh),\qquad [B]_G\longmapsto \left(\om_B\right)_{\mid \kh(x)}.
\]
On the other hand, \cite[Theorem B]{Rig} shows that the restriction mapping determines an isomorphism of posets:
\begin{equation}\label{thmB}
\V(\kh)\lra \V,\qquad \mu\longmapsto \mu_{\mid K(x)}.
\end{equation}
This ends the proof of our main theorem.

\begin{theorem}\label{main}
The mapping $\,\B/G\to \V$ determined by $[B]_G\mapsto \left(\om_B\right)_{\mid K(x)}$, is an isomorphism of posets.
\end{theorem}

In order to fit this theorem with Bengu\c{s}-Lasnier conjecture, we need only to identify $\B/G$ with some space $\D$ of diskoids. 

Let $\irr(\kh)$ be the set of all monic irreducible polynomials in $\kh[x]$.\e

\defn
Take $f\in\irr(\kh)$ and $\dta\in\qc$.
The {\bf diskoid} $\dfd$ centred at $f$ of radius $\dta$ is defined as
\[
\dfd=\{c\in\kb\mid\vb(f(c))\ge\dta^L\}.
\]
The {\bf pointed diskoid} $\dul(f,\dta)$ centred at $f$ of radius $\dta$ is defined as
\[
\dul(f,\dta)=\left(\dfd,\dta\right).
\]
We say that $\dfd$ is the {\bf underlying diskoid} of $\dul(f,\dta)$.

Let $\dqc$ be the set of all pointed diskoids centred at polynomials in $\irr(\kh)$ and having  radii in $\qc$.
Consider the following partial ordering on $\dqc$:
\[
\dul(f,\dta)\le \dul(g,\ep)\ \sii\ \dfd\supseteq D(g,\ep)\ \mbox{ and }\ \dta\le\ep
\]

Let $\dg,\dcut\subseteq \dqc$
be the subsets of all diskoids having radii in $\g$ and $\cuts$, respectively.

Let $\di$ be the quotient set of the set of nests of diskoids in $\dg$ having an empty intersection, under the equivalence relation identifying  mutually cofinal nests.

Let $\mathcal{D}=(D_i)_{i\in I}$ be a nest of diskoids, and $\dul\in\dqc$. Then, we define
\begin{equation}\label{partialDi}
\dul\le \mathcal{D}\ \sii\ \op{under}(\dul)\le D_i\ \mbox{ for some }i\in I, 
\end{equation}
where $\op{under}(\dul)$ is the underlying diskoid of $\dul$.
Clearly, this ordering is compatible with the equivalence  of nests of diskoids. 
 
Let $\D=\di\sqcup \dqc$ equipped with the partial ordering determined by the partial ordering of $\dqc$ and the relationship (\ref{partialDi}).

The following observation follows easily from \cite[Lemma 6.9]{Andrei}.
 
\begin{lemma}\label{B-D}
The following mapping is an isomorphism of posets:
\[
\B/G\lra \D,\qquad [B]_G\longmapsto \bigcup_{\sg\in G}\sg(B)
\]
\end{lemma}

Moreover, this isomorphism induces isomorphisms:
\[
\bi/G\iso \di,\quad \bg/G\iso\dg,\quad \bcut/G\iso\dcut,\quad \bqc/G\iso\dqc.
\] 

\subsection{Root node and leaves of $\V$}\label{subsecMM}

Clearly, the restriction mappings preserve the ordering of valuations. Hence, the restriction of the valuation $\omi$, defined in (\ref{omi}) is the root node of $\V$, and the restriction of the leaves of $\V(\kb)$ are the leaves of $\V$.

The set of leaves of $\V(\kb)$ was described in Section \ref{subsubsecMM}. Thus, $\V$ has a set of infinite leaves, parametrized by $\bi/G$, and a set of finite leaves, obtained as the restrictions to $K(x)$ of the valuations $\om_{a,\infty^-}$ for $a\in\kb$.

Since $\om_{a,\infty^-}$ and $\om_{b,\infty^-}$ have the same restriction to $\kh(x)$ if and only if $a$ and $b$ are $G$-conjugate, we see that the finite leaves of $\V(\kh)$ are
parametrized by the set $\irr(\kh)$. By the isomorphism (\ref{thmB}), the set of finite leaves of $\V$ is in bijection with $\irr(\kh)$ too.

\section{Valuations on $\kx$ with nontrivial support}\label{secSupport}

\subsection{Valuations on $\kx$}\label{secValsKx}
Let $\La$ be an ordered abelian group. A $\La$-valued valuation on the polynomial ring $\kx$ is a mapping $\mu\colon \kx\to \La\infty$, satisfying the following conditions:\e

(0) \ $\mu(1)=0$, \ $\mu(0)=\infty$,

(1) \ $\mu(fg)=\mu(f)+\mu(g),\quad\forall\,f,g\in \kx$,

(2) \ $\mu(f+g)\ge\min\{\mu(f),\mu(g)\},\quad\forall\,f,g\in \kx$.\e

The {\bf support} of $\mu$ is the prime ideal $\p=\supp(\mu)=\mu^{-1}(\infty)\in\op{Spec}(\kx)$.

The {\bf value group} of $\mu$ is the subgroup $\gm\subseteq \La$ generated by $\mu\left(\kx\setminus\p\right)$.

The valuation $\mu$ induces in a natural way a valuation $\mub$ on the field of fractions $L=\op{Frac}\left(\kx/\p\right)$; that is, $L=K(x)$ if $\p=0$, or $L=\kx/(f)$ if $\p=f\kx$ for some irreducible $f\in\kx$. 

The {\bf residue field} of $\mu$ is, by definition, the residue field $L\mub$.

Thus, the valuations on $\kx$ with trivial support may be identified with valuations on $K(x)$.
The valuations with nontrivial support extending $v$ are commensurable over $v$ and have a finite residual extension $L\mub/Kv$. In particular, these valuations are equivalent to a unique $\g$-valued valuation, by the arguments used in Section \ref{subsecComm}.

The valuations  with nontrivial support are useful to describe extensions of $v$ to finite field extensions of $K$. More precisely, every valuation $\mu$ with nontrivial support $f\kx$ can be identified with an extension of $v$ to the field $L=\kx/(f)$, just by considering the composition
\[
\mu\colon \kx\lra L\stackrel{\mub}\lra \g\infty.
\]

These extensions of $v$ to a given finite simple extension $L/K$, determined by an irreducible $f\in\kx$, are in 1-1 correspondence with the irreducible factors of $f$ in $\kh[x]$. For the discrete rank-one case this fact goes back to Hensel. For an arbitrary valued field $(K,v)$ it can be deduced from the techniques of \cite[Section 17]{endler}. A concrete proof may be found in \cite[Section 3]{NN}. 

For each irreducible factor $F\in\irr(\kh)$ of $f$, we may consider the following valuation $v_F$ on $\kx$ with support $f\kx$:
\[
v_F(g)=\vb(g(a))\quad\mbox{ for all }g\in\kx,
\] 
where $a\in\kb$ is any choice of a root of $F$ in $\kb$. By the henselian property, this valuation does not depend on the choice of $a$.  

As a consequence of the above mentioned 1-1 correspondence, the mapping $F\mapsto v_F$ determines a bijection between $\irr(\kh)$ and the set of all $\g$-valued valuations on $\kx$ with nontrivial support.

\subsection{Tree of equivalence classes of valuations on $\kx$}
Let $\T$ be the set of equivalence classes of valuations on $\kx$ whose restriction to $K$ is equivalent to $v$. This set has the structure of a  tree too, fully described in \cite[Section 7]{VT}. 

As a set, we have
\[
\T=\V\sqcup\{v_F\mid F\in \irr(\kh)\}.
\] 

What is the relative position of these added valuations with nontrivial support?

Recall that the finite leaves of $\V(\kb)$ are the valuations 
$\om_{a,\infty^-}$, for $a$ running on $\kb$.

For every $F\in\irr(\kh)$, the finite leaves $\om_{a,\infty^-}$, for $a$ running on the roots of $F$, have the same restriction to $\kh[x]$ (hence to $\kx$). This restriction satisfies
\[
\left(\om_{a,\infty^-}\right)_{\mid \kx}<v_F
\] 
and $v_F$ is the immediate successor of $\left(\om_{a,\infty^-}\right)_{\mid \kx}$ in $\T$.

Therefore, the finite leaves in $\V$ (parametrized by $\irr(\kh)$) cease to be leaves in $\T$, because each one admits an immediate successor in $\T$. 
Also, the valuations with nontrivial support become the set of finite leaves of $\T$.  

Finally, the tree $\T$ admits a completely analogous geometric parametrization as that of Theorem \ref{main}.

Let us briefly describe the parametrization of $\T(\kb)$. The parametrization of $\T$ follows exactly in the same way as before, just by considering $G$-action and applying \cite[Theorem B]{Rig}, which includes valuation with nontrivial support.  

The set $\bi$ remains untouched.
We enlarge the set $\bqc$, just by admitting radii in $\qc\infty$. We agree that
\[
B(a,\infty)=B(a,\infty^-)=\{a\}\quad\mbox{ for all }a\in \kb.
\]
The pointed balls $\bul(a,\infty)=(B(a,\infty),\infty)$ determine valuations $\om_{a,\infty}$ with nontrivial support $(x-a)\kb[x]$, by using exactly the same formula as in (\ref{defom}),  letting the symbol $x_\infty$ to be equal to $\infty$. That is,  
\[
\om_{a,\infty}\left(\sum\nolimits_{0\le n}a_n(x-a)^n\right) = \vb(a_0).
\]
The restriction of $\om_{a,\infty}$ to $\kx$ is the valuation $v_F$, for $F$ the minimal polynomial of $a$ over $\khx$.  

Denote by $\tilde{\B}_{\op{qcut}}$ this enlargement of $\bqc$. Then, the set
\[
\tilde{\B}:=\bi\sqcup \tilde{\B}_{\op{qcut}}
\]
parametrizes $\T(\kb)$ as indicated in Theorem \ref{mainB}.

\section{Comparison with approximation types}\label{secApprT}
The concept of approximation type was explored by Kuhlmann in \cite{KuhlAT}.

Let $(K,v)$ be a valued field with value group $\gv=vK$. Let $\mathfrak B=\mathfrak B(K,v)$ be the set of all balls (open or closed) in $K$, with radii in $\gv$. Let $\mathcal N$ be the set of all nests of balls in $K$. For $\mathcal B=(B_i)_{i\in I}\in\mathcal N$ we define
\[
\overline{\mathcal B}:=\{B\in \mathfrak B\mid B_i\subseteq B \ \mbox{ for some }i\in I\}.
\]

An {\bf approximation type} over $(K,v)$ is either the empty set or a nest of balls $\ap$ such that $\ap=\overline{\ap}$. In other words, a non-empty approximation type is a nest of balls containing every open or closed ball which contains some ball $B\in \ap$.

\begin{lemma}\label{AB}
 For every nest of balls $\mathcal B\in\mathcal N$, there exists a uniquely determined approximation type $\ap$ such that $\ap\sim \mathcal B$. 
\end{lemma}

\begin{proof}
	Let $\mathcal B$ be a nest of balls and set $\ap=\overline{\mathcal B}$. Clearly, $\ap=\overline{\ap}$, hence it is an approximation type. By definition, $\ap$ and $\mathcal B$ are cofinal in each other, so that $\ap\sim\bb$. Finally, two approximation types in $\mathcal N$ are mutually cofinal if and only if they are equal. Hence, $\textbf A$ is uniquely determined.
\end{proof}\e

The {\bf support} of $\ap$ is the following initial segment of $\gv$:
\[
\supp(\ap):=\left\{\ga\in\gv\mid \ap\ \mbox{ contains a closed ball of radius  }\ga\right\}. 
\]
For all $\ga\in\supp(\ap)$ there is a unique closed ball in $\ap$ of radius $\ga$. We denote this closed ball by $\ap_\ga$, and its corresponding open ball by $\ap_\ga^{\op{o}}$. 

Let $\mathcal A$ be the set of all approximation types over $(K,v)$. On this set, we consider the partial ordering determined by ascendent inclusion as sets of balls:
\[
\ap\le \ap'\ \sii\ \ap\subseteq \ap'.
\] 

\defn
We write $\bigcap\ap$ to indicate the intersection of all balls in $\ap$.
\begin{itemize}
\item $\ap$ is {\bf immediate}\, if $\ \ap\ne\emptyset$\, and $\ \bigcap\ap=\emptyset$.
\item  $\ap$ is {\bf value-extending}\, if either $\ap=\emptyset$, or 
\[
\bigcap\ap\ne\emptyset\quad \mbox{ and }\quad\ap_\ga^{\op{o}}\in\ap\ \mbox{ for all }\ga\in\supp(\ap).
\]
\item $\ap$ is {\bf residue-extending}\, if  there exists  $\ga\in\supp(\ap)$ such that $\ap_\ga^{\op{o}}\not\in\ap$.
\end{itemize}
In the latter case, we have $\ga=\max(\supp(\ap))$ and $\ \ap=\overline{\ap_\ga}$.\e

Let $\mathbb V$ be  the set of equivalence classes of valuations on $K(x)$ whose restriction to $K$ is equivalent to $v$.
For each valuation $\nu\in \mathbb V$, Kuhlmann defines the approximation type of $x$ over $(K,v)$ as follows
\[
{\rm appr}_\nu(x,K)=\{B\cap K\mid B\in\mathfrak B(K(x),\nu)\ \mbox{ centred on $x$, with radius in }\gv\}. 
\]
One can show that the above definition is compatible with equivalence of valuations on $K(x)$. The next result follows from \cite[Theorems 1.2, 1.3]{KuhlAT}.
\begin{theorem}\label{mainK}
	The map $\mathbb V\to \mathcal A$ given by $\nu\mapsto {\rm app}_\nu(x,K)$ is surjective.
	
	If $K$ is algebraically closed, then this mapping is bijective. 
\end{theorem}

Suppose from now on that $K$ is algebraically closed. Write $\g:=\gv$, which is a divisible group.

In this case, we have two different geometric parametrizations of $\V$. The mapping $\B\to\V$ from Theorem \ref{mainB} and the mapping $\V\to\mathcal{A}$ from Theorem \ref{mainK}. It is easy to check that the latter map strictly preserves the ordering, so that it is an isomorphism of posets too. Our aim in this section is to explicitly describe the isomorphism
\begin{equation}\label{composition}
\B\lra\V\lra \mathcal{A}
\end{equation}
obtained by composition of the two geometric parametrizations.

Recall the decomposition $\B=\bi\sqcup\bg\sqcup\bcut$ described in Section \ref{secKb}.

\begin{theorem}\label{comparison}
The isomorphism of (\ref{composition}) restricts to isomorphisms: 
\[
\bi\iso\mathcal{A}_\infty,\qquad \bg\iso\mathcal{A}_{\op{re}},\qquad \bcut\iso\mathcal{A}_{\op{ve}},
\]
where $\mathcal{A}_\infty$, $\mathcal{A}_{\op{re}}$ and $\mathcal{A}_{\op{ve}}$ are the subsets of immediate, residue-extending and value-extending approximation types.
More explicitly, the isomorphism acts as follows:
\[
\ars{1.3}
\begin{array}{lll}
\bi\lra \mathcal{A}_\infty,&\qquad [\bb]\mapsto \overline{\bb},&\\
\bg\lra \mathcal{A}_{\op{re}},&\qquad B(a,\ga)\mapsto \overline{B(a,\ga)},&\ \mbox{ for all }a\in K,\ \ga\in\g,\\
\bcut\lra \mathcal{A}_{\op{ve}},&\qquad B(a,\dta)\mapsto \overline{\left\{B^{\op{o}}(a,\ga)\mid \ga\in\dta^L\right\}},&\ \mbox{ for all }a\in K,\ \dta\in\cuts.
\end{array}
\]
\end{theorem}

\begin{proof}
Take any $[\bb]\in\bi$; that is, the class of a nest of closed balls $\bb=(B_i)_{i\in I}$ such that $\bigcap B_i=\emptyset$. Write $B_i=B(a_i,\gamma_i)$ for all $i\in I$, and let $S\subseteq \g$ be the initial segment generated by the subset $\{\ga_i\mid i\in I\}$. The set $S$ has no last element.

Let $\nu=\om_\bb$. By the remarks in Section \ref{subsecVAKb} we deduce that for every $b\in K$, there exists $i\in I$ such that
\[
\nu(x-b)= \min\{v(b-a_i),\gamma_i\}\leq \gamma_i.
\]
Hence, if $\gamma>S$, then $B(x,\ga)\cap K=\emptyset$. On the other hand, if $\ga<\ga_i$ for some $i\in I$, then $B^{\op{o}}(x,\ga)\cap K\supseteq B_i$. Indeed,
\[
b\in B_i \imp v(b-a_i)\ge \ga_i \imp \nu(x-b)\ge\om_{a_i,\ga_i}(x-b)=\min\{v(b-a_i),\ga_i\}>\ga.
\] 
Therefore, the image of $[\bb]$ in $\mathcal{A}$ is the approximation type
\[
\ap=\overline{\left\{B^{\op{o}}(x,\ga)\cap K\mid \ga\in S\right\}}.
\]
By Lemma \ref{AB}, in order to check that $\ap=\overline{\bb}$, we need only to show that $\ap\sim\bb$. We have already seen that every $B\in\ap$ contains some $B_i$. Let us show the converse statement; more precisely, let us show that every $B_i$ contains $B^{\op{o}}(x,\ga)\cap K$ for all $\ga\in S$ such that $\ga>\ga_i$. 

Indeed, take any $b\in B^{\op{o}}(x,\ga)\cap K$. Take $i<j\in I$ sufficiently large so that 
\[
\ga_i\le\nu(x-a_i)=\min\{v(a_j-a_i),\ga_j\},\qquad 
\ga<\nu(x-b)=\min\{v(a_j-b),\ga_j\}.
\]
The condition $v(b-a_i)<\ga_i$ leads to a contradiction:
\[
\ga_i>v(b-a_i)\ge\min\{v(b-a_j),v(a_j-a_i)\}\ge\min\{\ga,\ga_i\}=\ga_i.
\]
Hence, 
$B^{\op{o}}(x,\ga)\cap K\subseteq B_i$.

This proves that $\ap\sim\bb$. In particular, $\bigcap A=\bigcap_{i\in I}B_i=\emptyset$, so that $\ap$ is an immediate approximation type.
	
Now, take $B(a,\ga)\in\bg$ and let $\nu=\om_{a,\ga}$. Since $\nu(x-b)\le\ga$ for all $b\in K$, we have $B(x,\al)\cap K=\emptyset$ for all $\al>\ga$. On the other hand, for $\al=\ga$ we clearly have
\[
B(x,\ga)\cap K=B(a,\ga),\qquad B^{\op{o}}(x,\ga)\cap K=\emptyset.
\]
This proves that the image of $B(a,\ga)$ in $\mathcal{A}$ is the approximation type $\ap=\overline{B(a,\ga)}$, which is obviously residue-extending. 

Finally, take $B(a,\dta)\in\bcut$ for some $a\in K$ and $\dta\in\cuts$. Let $\nu=\om_{a,\dta}$. As in the previous cases, if $\ga\in\g$ satisfies $\ga>\dta^L$, then  $B(x,\ga)\cap K=\emptyset$. On the other hand, if $	\ga\in\dta^L$, then $B^{\op{o}}(x,\dta)\cap K=B^{\op{o}}(a,\dta)$. Therefore,  the image of $B(a,\dta)$ in $\mathcal{A}$ is the approximation type
\[
\overline{\left\{B^{\op{o}}(a,\ga)\mid \ga\in\dta^L\right\}},
\]
which is obviously value-extending. 
\end{proof}\e

Let us illustrate with some examples, the fact that the isomorphisms in (\ref{composition}) preserve the ordering.

The empty approximation type belongs to $\mathcal{A}_{\op{ve}}$. It corresponds to the maximal ball $B(a,-\infty)=K\in\bcut$, which is independent of $a\in K$. They both correspond to the root node $\omi$ of $\V$.

The leaves of $\V$ correspond to the minimal balls $B(a,\infty^-)=\{a\}\in\bcut$ for all $a\in K$. The corresponding  maximal approximation types in $\mathcal{A}_{\op{ve}}$ are  those of the form 
\[
\overline{\left\{B^{\op{o}}(a,\ga)\mid \ga\in\g\right\}}=\left\{B\in\mathfrak{B}\mid a\in B\right\}.
\]

Finally, take arbitrary $a\in K$, $\ga\in \g$. Consider the quasi-cuts $\ga^-<\ga<\ga^+$, and the valuations
$\om_{a,\ga^-}<\om_{a,\ga}<\om_{a,\ga^+}$. 
The  corresponding approximation types are
\[
\overline{\left\{B^{\op{o}}(a,\al)\mid \al<\ga\right\}}\subsetneq \overline{B(a,\ga)}\subsetneq \overline{B^{\op{o}}(a,\ga)}.
\]

\end{document}